\documentclass[11pt]{amsart}

\usepackage{amsmath,amstext,amsthm,amssymb}

\begin{document}

\newcommand{\N}{\mathbb{N}}
\newcommand{\Z}{\mathbb{Z}}
\newcommand{\R}{\mathbb{R}}
\newcommand{\C}{\mathbb{C}}
\newcommand{\Q}{\mathbb{Q}}

\newcommand{\dlat}{\mathfrak{d}}
\newcommand{\otau}{\overline{\tau}}

\newcommand{\diag}{\mathrm{diag}}
\newcommand{\vol}{\mathrm{vol}}
\newcommand{\lin}{\mathrm{lin}}
\newcommand{\inter}{\mathrm{int}}
\newcommand{\cl}{\mathrm{cl}}
\newcommand{\bd}{\mathrm{bd}}
\newcommand{\relbd}{\mathrm{relbd}}
\newcommand{\relint}{\mathrm{relint}}
\newcommand{\conv}{\mathrm{conv}}
\newcommand{\ogamma}{\overline{\gamma}}
\newcommand{\raninv}{\tau}

\newcommand{\pos}{\mathrm{pos}}
\newcommand{\aff}{\mathrm{aff}}
\newcommand{\co}{\mathrm{c}}
\newcommand{\coh}{\mathrm{h}}

\newcommand{\e}{{\mathrm{e}}}
\newcommand{\ic}{{\mathrm{i}}}
\newcommand{\ov}{\overline}
\newcommand{\dist}{\mathrm{dist}}
\newcommand{\lwidth}{\mathrm{w}}
\newcommand{\norcone}{\mathcal{N}}
\newcommand{\admi}{\mathcal{A}}
\newcommand{\denl}{\delta_\Lat}
\newcommand{\K}{\mathcal{K}}
\newcommand{\cK}{\mathcal{K}_0^n}
\newcommand{\Lat}{\mathcal{L}}

\newcommand{\dual}[1]{{#1^\mathfrak{*}}}

\newcommand{\norm}[2]{\Vert #1\Vert_#2}
\newcommand{\enorm}[1]{\Vert #1\Vert}
\newcommand{\onorm}[1]{\Vert #1\Vert_1}
\newcommand{\Cnorm}[1]{\Vert #1\Vert_\infty}
\newcommand{\grass}{\mathrm{G}}
\newcommand{\four}[1]{\widehat{#1}}
\newcommand{\facet}{\mathcal{F}}
\newcommand{\pack}{\mathcal{P}}
\newcommand{\packl}{\mathcal{P}_\Lat}
\newcommand{\bfill}{\mathfrak{s}}
\newcommand{\DV}{\mathfrak{V}}
\newcommand{\Lamp}{\Lambda^\circ}
\newcommand{\trans}{\intercal}

\newcommand{\disuni}{\mathbin{\setbox0\hbox{$\bigcup$}\rlap{\copy0}\raise.3%
  \ht0\hbox to \wd0{\hfil$\cdot$\hfil}}}
\theoremstyle{plain}
\newtheorem{theorem}{Theorem}[section]
\newtheorem{corollary}{Corollary}[section]
\newtheorem{proposition}{Proposition}[section]
\newtheorem{remark}{Remark}[section]
\newtheorem {definition}{Definition}[section]
\newtheorem {example}{Example}[section]
\newtheorem {algorithm}{Algorithm}[section]
\newtheorem {problem}{Problem}[section]
\newtheorem {examples}{Examples}[section]
\newtheorem {lemma}{Lemma}[section]
\newtheorem {observation}{Observation}[section]
\numberwithin{equation}{section}

\title{Free planes in lattice sphere packings}
\author{Martin Henk}
\address{Universit\"at Magdeburg, Institut f\"ur Algebra und Geometrie,
  Universit\"ats\-platz 2, D-39106 Magdeburg, Germany}
\email{henk@math.uni-magdeburg.de}

\begin{abstract}
  We show that for every lattice packing of $n$-dimensional spheres 
  there exists an $(n/\log_2(n))$-dimensional affine plane which does not
  meet any of the spheres in their interior, provided $n$ is large
  enough. Such an affine plane is called a free plane and our result
  improves on former bounds.   
 \end{abstract}
\maketitle
\section{Introduction}
The main purpose of this note is to show that the complement of every
$n$-dimensional lattice sphere packing contains an $n/\log_2(n)$ affine
plane, if the dimension $n$ is large enough. In order to show the
relations of this result to other problems in packing theory we 
give a more general introduction. 

Let $\cK$ be the set of all $0$-symmetric convex bodies in the
$n$-dimensional Euclidean space $\R^n$. We always assume that
$\inter(K)\ne\emptyset$ for $K\in \cK$, i.e.,  $K$ has non-empty
interior, and we denote the volume of $K$ ($n$-dimensional Lebesgue
measure) by $\vol(K)$. A lattice $\Lambda\subset\R^n$ is the image  
of the standard lattice $\Z^n$ under a regular linear map $B$, i.e.,
$\Lambda=B\Z^n$, $B\in\R^{n\times n}$ with $\det B\ne 0$. $|\det B|$
is called the determinant of $\Lambda$ and it is denoted by
$\det\Lambda$. $\Lat^n$ denotes the set of all lattices in $\R^n$.

For a given $K\in\cK$, a lattice $\Lambda\in\Lat^n$ is called a packing
lattice of $K$ if for $b_1\ne b_2\in\Lambda$ the translates $b_1+K$
and  $b_2+K$ do not overlap, i.e.,
$\inter(b_1+K)\cap\inter(b_2+K)=\emptyset$. The density of a densest  lattice
packing of $K$, denoted by $\delta(K)$, is the maximum 
proportion of space that can be occupied by an arrangement of the
 $\Lambda+K$, where $\Lambda$ is a packing lattice. Thus 
\begin{equation*}
  \delta(K) = \max\left\{\frac{\vol(K)}{\det\Lambda}:
    \Lambda\hbox{ packing lattice of } K\right\}.
\end{equation*}
For a
thorough treatment of  packings of convex bodies and lattices we refer
to  \cite{GruLek:geonum} and \cite{ErdGruHam:lattpoint}. 

For a given lattice $\Lambda$ and a body $K\in\cK$ the maximum
dilation factor $\lambda$ with the property that $\Lambda$ becomes a packing lattice
of $\lambda\,K$ is called the packing radius of $\Lambda$ with respect to
$K$ and is denoted by $\lambda(K,\Lambda)$, i.e,
\begin{equation*}
  \lambda(K,\Lambda)=\max\left\{\lambda\in\R_{>0} : \Lambda \hbox{ is
      a packing lattice of }\lambda\,K\right\}.
\end{equation*} 
In particular, $\Lambda/\lambda(K,\Lambda)$ is the 
``smallest'' dilation of $\Lambda$ which is a packing lattice of
$K$. Hence we can write 
\begin{equation*}
 \delta(K)=\vol(K)\,
\max\left\{\frac{\lambda(K,\Lambda)^n}{\det\Lambda}: \Lambda\in\Lat^n \right\},
\end{equation*}
and thus 
\begin{equation}
  \lambda(K,\Lambda)^n \leq \delta(K)\,\frac{\det\Lambda}{\vol(K)}.
\label{eq:bound_on_packing_radius}
\end{equation}
As a counterpart to the packing radius we have the covering radius 
$\mu(K,\Lambda)$ of $\Lambda$ with respect to $K$, defined by 
\begin{equation*}
  \mu(K,\Lambda)=\min\left\{\mu\in\R_{>0} : \Lambda + \mu\,K=\R^n\right\}.
\end{equation*}
One can also say, that $\mu(K,\Lambda)$ is the smallest positive
number 
$\mu$ such that  $\Lambda$ is a covering lattice of
$\mu\,K$. If $\theta(K)$ denotes the density of a thinnest lattice
covering of $K$,  we get analogously to
\eqref{eq:bound_on_packing_radius} 
\begin{equation}
 \mu(K,\Lambda)^n \geq \theta(K)\,\frac{\det\Lambda}{\vol(K)}.
\label{eq:bound_on_covering_radius}
\end{equation}

Here we are interested in lattices which are simultaneously
good packing and covering lattices. Therefore we define  
\begin{equation*}
   \gamma(K)=\inf\left\{ \mu(K,\Lambda) : \Lambda \hbox{ is a packing
       lattice of }K\right\}.
\end{equation*} 
{\sc Zong} baptised  this number {\em the simultaneous lattice
  packing and covering constant of $K$} (cf.~\cite{zong:bams}). By definition of the packing
radius of $K$ we may also write 
\begin{equation*}
   \gamma(K)=\inf\left\{ \frac{\mu(K,\Lambda)}{\lambda(K,\Lambda)} :
     \Lambda \in\Lat^n\right\}
\end{equation*} 
and by \eqref{eq:bound_on_packing_radius} and
\eqref{eq:bound_on_covering_radius} we get 
\begin{equation}
 \gamma(K)\geq \left(\frac{\theta(K)}{\delta(K)}\right)^{1/n}.
\label{eq:lower_bound_gamma}
\end{equation}
A remarkable upper bound on $\gamma(K)$ is due to {\sc
  Butler} \cite{But:paccov}, who showed  
\begin{equation*}
     \gamma(K) \leq 2 +o(1),
\end{equation*}
as $n$ tends to infinity. By \eqref{eq:lower_bound_gamma} it gives
asymptotically a lower bound on $\delta(K)$ which is of the same order
of magnitude as the classical Minkowski-Hlawka bound. 

$\gamma(K)$ is a quite interesting functional. If we were able to  prove, for
instance 
for the unit ball $B^n$, 
the existence of an $\epsilon >0$ such that
$\gamma(B^n)\leq 2-\epsilon$ then we would get
by \eqref{eq:lower_bound_gamma} an improvement on the
best known lower bound on $\delta(B^n)$ which is of order
$2^{-n}$ \cite{Ball:sphere_packings}.

On the other hand, if we could verify $\gamma(K)\geq 2$ then we know
that even for a densest packing lattice $\Lambda_K$ of $K$, i.e,
$\delta(K)=\vol(K)/\det\Lambda_K$,  there exists a point
$x\in\R^n$ with  $(x+\inter (K))\cap
(\Lambda_K+\inter(K))=\emptyset$.  Thus $\Lambda_K\cup(x+\Lambda)$
is a packing set, which means that the bodies
$(\Lambda_K\cup(x+\Lambda_K))+K$ do not overlap. 
Since the space is
occupied by $(\Lambda_K\cup(x+\Lambda_K))+K$ two times as good as by
$\Lambda_K+K$ we get  
\begin{equation*}
   \delta_T(K) \geq 2\delta(K),  
\end{equation*}  
where $\delta_T(K)$ denotes the density of a densest arbitrary
(without the restriction to lattices) packing of $K$. So far a body
$K\in\cK$ with $\delta_T(K)>\delta(K)$ is not known. 
For results on $\gamma(K)$ in the planar and three dimensional case we
refer to \cite{zong:sim_2} and \cite{zong:sim_3}.
 
Next we are interested in a generalisation of $\gamma(K)$. To this end
we consider the covering radii introduced
by {\sc Kannan and Lov\'asz} \cite{kannan_lovasz:covering_minima}. For
$1\leq i\leq n$, 
\begin{equation*}
\begin{split}
  \mu_i(K,\Lambda) = \min\big\{\mu >0  : &\, \Lambda+\mu K \text{ meets
      every $(n-i)$-dimensional} \\ 
      & \text{affine subspace of $\R^n$}\big\}
\end{split}
\end{equation*}
is called the $i$-th covering minimum. With the help of theses
functionals we define for $1\leq i\leq n$ 
 \begin{equation*}
   \gamma_i(K)=\inf\left\{ \frac{\mu_i(K,\Lambda)}{\lambda(K,\Lambda)} : \Lambda \in\Lat^n\right\}.
\end{equation*} 
Of course, $\gamma_n(K)=\gamma(K)$, and by definition we have 

\vspace{0.5cm}
\noindent
\begin{minipage}{3.5cm}
\begin{center}
\begin{equation}
\gamma_i(K)\geq 1 \, \Leftrightarrow
\label{eq:free_planes}
\end{equation}
\end{center}
\end{minipage}
\begin{minipage}{8.9cm}
For every lattice packing $\Lambda$ of $K$ there exists an
$(n-i)$-dimensional affine plane which does not intersect any of the translates $\Lambda+K$ in their interior.
\end{minipage}
\vspace{0.5cm}


A plane which does not intersect any of the translates of the packing is called a free plane. In general, we can not expect to find 
a free plane for an arbitrary $K\in\cK$ and $\Lambda\in\Lat^n$. However, it was shown by  Heppes
 \cite{Heppes:free_planes} 
that in every 3-dimensional lattice sphere packing one can find a 
cylinder of infinite length which  does not intersect any of the
spheres of this packing. In other word, one can always find an
$1$-dimensional free plane. Later this was generalised to any dimension $\geq 3$ by Horv\'ath and Ry{\v s}kov \cite{Horvath_Ryshkov:free_planes}.
Concerning the maximum dimension of free planes in lattice sphere
packings it was recently shown 
that  it is at least $4$ ($n$ large) and at most $c\,n$ for an absolute
constant $c<1$ \cite{henk_ziegler_zong:free}. For a detailed treatment of the history of free planes
in sphere packings 
and related results we refer to the recent survey of {\sc Zong} \cite{zong:bams}.

Here we improve on the upper bound on free planes in the following way
\begin{theorem} For every lattice packing of $n$-dimensional spheres, $n$ large, 
  there exists an affine plane of dimension $\geq n/\log_2(n)$ which
  does not intersect any of the spheres of the packing in their
  interior. 
\label{thm:main}
\end{theorem}

By \eqref{eq:free_planes} the theorem will be an immediate consequence
of the next lemma giving a lower bound on $\gamma_i(K)$. 
To this end, for $1\leq i\leq n$ and $K\in\cK$ let 
\begin{equation*}
   \rho_i(K)=\sup_{A\in\mathrm{GL}(n,\R)}\, \inf \left\{\frac{\vol(AK)^{1/n}}{\vol_i(\,(AK)|L\,)^{1/i}}:L\in \mathrm{G}(i,n) \right\}.
\end{equation*}
Here $\mathrm{GL}(n,\R)$ denotes the set of all regular $(n\times n)$-matrices, $\mathrm{G}(i,n)$ denotes the family of all $i$-dimensional linear subspaces of $\R^n$, and for a set $S$ the orthogonal projection onto a linear subspace $L$ is denoted by $S|L$. The $i$-dimensional volume of an $i$-dimensional  $S\subset\R^n$ is denoted by $\vol_i(S)$.
\begin{lemma} Let $K\in\cK$. Then 
\begin{equation*}
 \gamma_i(K) \geq \delta(K)^{-\frac{1}{n}}\cdot\rho_{i}(K)\cdot   
 \,
 \max
 \left\{ \left(\frac{1}{\sqrt{i+1}}\right)^\frac{n-i}{i}, 
         \frac{1}{\sqrt{n-i+1}}\right\}.
\end{equation*}   
In particular, for the sphere and $n$ large we have 
\begin{equation*}
 \gamma_i(B^n) > \frac{3}{2}\,\sqrt{\frac{i}{n}}\,  \max
 \left\{ \left(\frac{1}{\sqrt{i+1}}\right)^\frac{n-i}{i}, 
         \frac{1}{\sqrt{n-i+1}}\right\}.
\end{equation*}
\label{lem:main}
\end{lemma}

\vspace{0.3cm}
\noindent
Substituting $i=n-n/\log_2(n)$ in the last formula leads to 
\begin{equation*}
  \gamma_{n-\frac{n}{\log_2(n)}}(B^n) >
  \frac{3}{2}\,\sqrt{1-\frac{1}{\log_2(n)}}\,
  \sqrt{\frac{\log_2(n)}{n\,\log_2(n)-n+\log_2(n)}}^{\,\,\frac{1}{\log_2(n)-1}}.                    
\end{equation*}  
Since the last factor tends to $1/\sqrt{2}$ as $n$ approaches infinity,
Theorem \ref{thm:main} follows from the observation
\eqref{eq:free_planes}. 

Before giving the proof of the lemma we want to remark that in
every dimension there exists a lattice $\tilde{\Lambda}\in\Lat^n$ such
that $\mu_i(B^n,\tilde{\Lambda})/\lambda(B^n,\tilde{\Lambda})\leq \co\, \frac{i}{n}$,
where $\co>2$ is an absolute constant
(cf.~\cite{henk_ziegler_zong:free}). Hence, at least for the sphere
and $i=1$ Lemma \ref{lem:main} is best possible up to a constant.

\section{Proof}
First we have to fix some notations. A linear $i$-dimensional subspace
$L_i\in \mathrm{G}(i,n)$ is called a lattice plane (subspace) of
$\Lambda\in\Lat^n$ if $\dim(\Lambda\cap L_i)=i$. The set of all
$i$-dimensional lattice planes of a given lattice $\Lambda$ is denoted
by $\mathcal{L}(\Lambda,i)$. 
 For a linear subspace $L\subset\R^n$ we denote by $L^\perp$ the
orthogonal complement. We note that the orthogonal projection of
a lattice $\Lambda\in\mathcal{L}^n$ onto a $i$-dimensional linear
subspace $L$, say, is a lattice again, if and only if $L\in
\mathcal{L}(\Lambda,i)$.  The polar (dual) lattice of $\Lambda\in\Lat^n$ is
denoted by $\dual{\Lambda}$ and it is 
$\det\dual{\Lambda}=1/\det\Lambda$. There are some  simple and   
useful relations 
between a lattice $\Lambda$ and its dual 
$\dual{\Lambda}$, which we collect in the following statements   
\begin{equation}
\begin{split}
  {\rm i)}&\quad L \in\mathcal{L}(\Lambda,i) \Longleftrightarrow L^\perp \in
  \mathcal{L}(\dual{\Lambda},n-i), \\
  {\rm ii)}&\quad \dual{(\Lambda| L^\perp)}=\dual{\Lambda}\cap
  L^\perp,\quad  L \in\mathcal{L}(\Lambda,i), \\
  {\rm iii)}&\quad \det(\Lambda\cap L) =
  \det(\dual{\Lambda}\cap L^\perp)\cdot\det\Lambda,\quad  L \in\mathcal{L}(\Lambda,i).
\end{split} 
\label{eq:proj_sect}
\end{equation}

For $1\leq i\leq n$ let 
\begin{equation*}
       \raninv(n,i) = \sup_{\Lambda\in\mathcal{L}^n} \inf \left\{
       \frac{\det(\Lambda\cap L)^{1/i}}{\det\Lambda^{1/n}} : 
        L\in\mathcal{L}(\Lambda,i)\right\}.
\end{equation*}
 The functional $[\raninv(n,i)]^{2\,i}$ is
the so    
called $i$-th Rankin invariant (or generalised Hermite constant),  
introduced by Rankin \cite{Ran:general_hermite} (see also
\cite{mullender:rankin}).  For our
purposes, however, the above normalisation is more suitable. 
Among others Rankin proved two basic relations for these numbers, which read in our notation 
\begin{equation}
 \begin{split}
 {\rm i)}&\, \raninv(n,i)^i=\raninv(n,n-i)^{n-i},\\
 {\rm ii)}&\, \raninv(n,i)\leq \raninv(n,m)\raninv(m,i),\quad\, 1\leq i\leq
 m\leq n.
 \end{split}
\label{eq:rankin}
\end{equation}
In fact, the identity i) is a consequence of
\eqref{eq:proj_sect} iii), and the second inequality  follows immediately from the
definition. 

Of particular interest is the $\raninv(n,1)$ which
can be rewritten as  
\begin{equation}
     \raninv(n,1)=2\left(\frac{\delta(B^n)}{\kappa_n}\right)^{1/n},
\label{eq:constd1}
\end{equation}
where 
$\kappa_n=\pi^{n/2}/\Gamma(n/2+1)$ is the volume of the $n$-dimensional unit ball $B^n$.
It is well known that (cf.~e.g.~\cite[p.~53]{martinet:lattices})  
\begin{equation}
               \raninv(n,1)<\sqrt{n}.
\label{eq:bound_d1}
\end{equation}
As a simple consequence of  properties \eqref{eq:rankin} we get for
$1\leq i\leq n$ 
\begin{equation}
\raninv(n,i) \leq
     \min\left\{\sqrt{n-i+1},\,\,
     \left(\sqrt{i+1}\right)^\frac{n-i}{i}\right\}.
\label{eq:simple_rankin}
\end{equation}
In fact, let  $i\leq n-1$. From 
\eqref{eq:rankin} ii) and \eqref{eq:bound_d1}  we obtain 
\begin{equation*}
\begin{split}
  \raninv(n,i)&\leq \prod_{m=i}^{n-1} \raninv(m+1,m) = \prod_{m=i}^{n-1}
  \raninv(m+1,1)^\frac{1}{m} \leq 
  \prod_{m=i}^{n-1} \left(\sqrt{m+1}\right)^\frac{1}{m} \\ &\leq 
  \left(\sqrt{i+1}\right)^\frac{n-i}{i}. 
\end{split}
\label{eq:produ}
\end{equation*}
Next we 
apply this  bound  to $\raninv(n,n-i)$ and with \eqref{eq:rankin}
i) we find  
$$
  \raninv(n,i)= \raninv(n,n-i)^\frac{n-i}{i} \leq \sqrt{n-i+1}.
$$
and so we get \eqref{eq:simple_rankin}.
\begin{proof}[Proof of Lemma \ref{lem:main}] Since $\gamma_i(K)$ is
  invariant with respect to linear transformations we may assume that
  $K$ is in such a position that 
\begin{equation*}
  \sup_{A\in\mathrm{GL}(n,\R)}\, \inf_{L\in \mathrm{G}(n,i)} \frac{\vol(AK)^{1/n}}{\vol_i(\,(AK)|L\,)^{1/i}}
\end{equation*} 
is attained for the identity matrix. 
It was shown in
  \cite{kannan_lovasz:covering_minima} that the $i$-th covering
  minimum can also be described as 
$$
  \mu_i(K,\Lambda)=\max\left\{ \mu(K|L^\perp,\Lambda|L^\perp) : 
                        L\in\mathcal{L}(\Lambda,n-i)\right\}.
$$
Hence  we have $\mu_i(K,\Lambda)^i\,\vol_i(K|L^\perp) \geq \det(\Lambda|L^\perp)$
for all $L\in\mathcal{L}(\Lambda,n-i)$ and with \eqref{eq:proj_sect}
we obtain 
\begin{equation*}
 \begin{split}
    &\frac{\mu_i(K,\Lambda)}{\lambda(K,\Lambda)} \geq
  \frac{1}{\lambda(K,\Lambda)}\max\left\{
    \frac{\det(\Lambda|L^\perp)^{1/i}}{\vol_i(K|L^\perp)^{1/i}} : L\in\mathcal{L}(\Lambda,n-i)\right\}
  \\ 
\quad &= \frac{1}{\lambda(K,\Lambda)}\,
      \max\left\{ \frac{\det(\Lambda)^{1/n}}{\vol_i(K|L)^{1/i}}\frac{\det(\dual{\Lambda})^{1/n}}
                       { \det(\dual{\Lambda}\cap L)^{1/i}} :
                L\in\mathcal{L}(\dual{\Lambda},i)\right\} \\
\quad &= \frac{\det(\Lambda)^{1/n}}{\vol(K)^{1/n}\,\lambda(K,\Lambda)}\,
      \max\left\{ \frac{\vol(K)^{1/n}}{\vol_i(K|L)^{1/i}}\frac{\det(\dual{\Lambda})^{1/n}}
                       { \det(\dual{\Lambda}\cap L)^{1/i}} :
                L\in\mathcal{L}(\dual{\Lambda},i)\right\}
\\
\quad  &\geq \delta(K)^{-1/n}\,\rho_i(K)\,\tau{(n,i)}^{-1}.
 \end{split}
\label{eq:main_eq}
\end{equation*}
Together with \eqref{eq:constd1} we obtain Lemma \ref{lem:main} in the general case. 
Next we firstly observe that by the isoperimetric inequality $\rho_i(K)\leq \kappa_n^{1/n}/\kappa_i^{1/i}$ with equality if and only if $K$ is an ellipsoid. Hence for 
$K=B^n$ we have 
\begin{equation*}
  \rho_i(B^n)=\frac{\kappa_n^{1/n}}{\kappa_i^{1/i}}=\frac{\Gamma(n/2+1)^{1/n}}{\Gamma(i/2+1)^{1/i}}\geq \sqrt{\frac{i}{n}},
\end{equation*} 
where the last inequality follows form the fact that $\Gamma(m/2+1)^{1/m}/\sqrt{m}$ is monotonously decreasing which can be shown by elementary considerations. Finally, as an upper bound on $\delta(B^n)$ we use the  
asymptotic bound of Kabatianski
and  Leven{\v s}tein $\delta(B^n)\leq 2^{-0.599\,n+o(1)}$ (cf.~\cite[p.~137]{zong:sphere}).
\end{proof}

\section{Concluding remarks} 
We remark that 
already Rankin determined the value $\raninv(4,2)=(3/2)^{1/4}$. For further
bounds for special choices of $i$ and $n$ and the relation of Rankin's
invariants to minimal vectors of the $i$-th exterior power of a
lattice we refer to \cite{coulangeon:hermite} and the references
within. Thunder
\cite{thunder:gen_hermite} gave a  lower bound on a generalised version 
of $\raninv(n,i)$ in the context of  
algebraic extensions $F$ of $\Q$ of finite degree. In the case $F=\Q$ his bound becomes 
\begin{equation*}
 \raninv(n,i)\geq \left(2\,n\,\frac{\prod_{j=1}^i\frac{\zeta(n-i+j)}{(n-i+j)\kappa_{n-i+j}} }{\prod_{j=2}^i\frac{\zeta(j)}{j\,\kappa_j} }\right)^\frac{1}{i\,n},
\end{equation*} 
where $\zeta(k)=\sum_{i=1}^\infty i^{-k}$ denotes the $\zeta$-function. For $i$ proportional to $n$ and $n\to\infty$ this lower bound is of order $\sqrt{n}^{(n-i)/n}$.
However, even if this lower bound were the right order of magnitude for an upper bound on $\raninv(n,i)$ this would not lead to an improvement on the  lower bound on the maximum dimension of free planes in lattice sphere packings. 

Finally, we note that free planes are lattice phenomena. In the non-lattice case one can even find sphere packings such that the length of a segment contained in the complement is bounded by a constant depending only on the dimension (see \cite{BorTarods:segments_complement} and \cite{henk_zong:segments}). 

\vspace{1cm}
\noindent 
{\em Acknowledgement.} I would like to thank K\'aroly B\"or\"oczky, Jr. and Achill Sch\"urmann for valuable comments and discussions. 
\providecommand{\bysame}{\leavevmode\hbox to3em{\hrulefill}\thinspace}
\providecommand{\MR}{\relax\ifhmode\unskip\space\fi MR }
\providecommand{\MRhref}[2]{%
  \href{http://www.ams.org/mathscinet-getitem?mr=#1}{#2}
}
\providecommand{\href}[2]{#2}



\begin{thebibliography}{Zon02b}

\bibitem[Bal92]{Ball:sphere_packings}
K.~Ball, \emph{A lower bound for the optimal density of lattice packings},
  International Math. Research Notes \textbf{10} (1992), 217--221.

\bibitem[BT]{BorTarods:segments_complement}
K.~B{\"o}r{\"o}czky, Jr. and G.~Tardos, \emph{The longest segment in the
  complement of a packing}, Mathematika, to appear.

\bibitem[But72]{But:paccov}
G.J. Butler, \emph{Simultaneous packing and covering in euclidean space},
  Proc.~LOndon Math.~Soc.~(3) \textbf{25} (1972), 721--735.

\bibitem[Cou97]{coulangeon:hermite}
Renaud Coulangeon, \emph{Minimal vectors in the second exterior power of a
  lattice}, Journal of algebra \textbf{194} (1997), 467--476.

\bibitem[EGH89]{ErdGruHam:lattpoint}
P.~Erd{\H o}s, P.M. Gruber, and J.~Hammer, \emph{Lattice points}, Longman
  Scientific \& Technical, Harlow, Essex/Wiley, New York, 1989.

\bibitem[GL87]{GruLek:geonum}
P.M. Gruber and C.G. Lekkerkerker, \emph{Geometry of {N}umbers}, 2nd ed.,
  North-Holland, Amsterdam, 1987.

\bibitem[Hep61]{Heppes:free_planes}
A.~Heppes, \emph{Ein {S}atz \"uber gitterf\"ormiger {K}ugelpackungen},
  Ann.Univ. Sci. Budapest E\"otv\"os Sect. Math. (1960/61), 89--90.

\bibitem[HR75]{Horvath_Ryshkov:free_planes}
J.~Horv\'ath and S.S. Ry{\v s}kov, \emph{Estimation of the radius of a cylinder
  that can be imbedded in any lattice packing of $n$-dimensional unit balls},
  Mat. Zametki \textbf{17} (1975), 123--128.

\bibitem[HZ00]{henk_zong:segments}
M.~Henk and C.~Zong, \emph{Segments in ball packings}, Mathematika \textbf{47}
  (2000), 31--38.

\bibitem[HZZ02]{henk_ziegler_zong:free}
M.~Henk, G.M. Ziegler, and C.~Zong, \emph{On free planes in lattice ball
  packings}, Bull. London Math. Soc. \textbf{34} (2002), no.~3, 284--290.

\bibitem[KL88]{kannan_lovasz:covering_minima}
R.~Kannan and L.~Lov\'asz, \emph{Covering minima and lattice-point-free convex
  bodies}, Annals of Mathematics \textbf{128} (1988), 577--602.

\bibitem[Mar03]{martinet:lattices}
J.~Martinet, \emph{Perfect lattices in {E}uclidean spaces}, Springer, Berlin
  Heidelberg, 2003.

\bibitem[Mul64]{mullender:rankin}
P.~Mullender, \emph{Some remarks on a method of mordell in the geometry of
  numbers}, Acta Arithm. \textbf{9} (1964), 301--304.

\bibitem[Ran53]{Ran:general_hermite}
R.A. Rankin, \emph{On positive definite quadratic forms}, J. London Math. Soc.
  \textbf{28} (1953), 309--314.

\bibitem[Thu98]{thunder:gen_hermite}
J.~L. Thunder, \emph{Higher dimensional analogs of hermite's constants},
  Michigan Math. J. \textbf{45} (1998), no.~2, 301--314.

\bibitem[Zon99]{zong:sphere}
Chuanming Zong, \emph{Sphere packings}, Universitext, Springer-Verlag, New
  York, 1999.

\bibitem[Zon02a]{zong:bams}
\bysame, \emph{From deep holes to free planes}, Bull. Amer. Math. Soc. (N.S.)
  \textbf{39} (2002), no.~4, 533--555 (electronic).

\bibitem[Zon02b]{zong:sim_2}
\bysame, \emph{Simultaneous packing and covering in the {E}uclidean plane},
  Monatsh. Math. \textbf{134} (2002), no.~3, 247--255.

\bibitem[Zon03]{zong:sim_3}
\bysame, \emph{Simultaneous packing and covering in three-dimensional
  {E}uclidean space}, J. London Math. Soc. (2) \textbf{67} (2003), no.~1,
  29--40.

\end{thebibliography}
\end{document}